\documentclass[12pt]{article}
\usepackage[english]{babel}
\usepackage[top=3.0cm, left=3.0cm, bottom=3.0cm, right=3.0cm]{geometry}
\date{\ }
\usepackage{enumerate}

\usepackage{float}

\usepackage{epsfig,graphicx,graphics,latexsym,amssymb,amsfonts,amsmath,amscd, euscript}
\usepackage{changebar}
\usepackage{graphicx}

\newtheorem{definition}{{\bf Definition}}[section]
\newtheorem{theorem}[definition]{{\bf Theorem}}

\newtheorem{corollary}[definition]{{\bf Corollary}}

\newtheorem{lemma}[definition]{\noindent {\bf Lemma}}

\newtheorem{claim}[definition]{\noindent {\bf Claim}}

\def\Proof{{\parindent0pt {\bf Proof.\ }}}

\def\endproof{\hfill {\kern 6pt\penalty 500
   \raise -0pt\hbox{\vrule \vbox to5pt {\hrule width 5pt
  \vfill\hrule}\vrule}}}
\usepackage{enumerate}

\newcommand{\K}{\mathbb{K}}

\newcommand{\Z}{\mathbb{Z}}
\newcommand{\Q}{\mathbb{Q}}

\begin{document}

\title{\bf On a generalization of Kelly's combinatorial lemma}
\author{Aymen BEN AMIRA $^1$, Jamel DAMMAK $^1$, Hamza SI KADDOUR $^{2, {\ast}}$\\
\centerline{$^{1}$  {\small Department of Mathematics, Faculty of Sciences of Sfax, B.P. 802, 3018 Sfax, Tunisia}}\\
\centerline{$^{2}$  {\small ICJ, Department of Mathematics, University of Lyon, University Claude-Bernard  Lyon1,}}\\
\centerline{\small 43 Bd du 11 Novembre 1918, 69622 Villeurbanne Cedex, France}}

\maketitle
 \footnotetext{$^{\ast}$ Correspondence: sikaddour@univ-lyon1.fr}

\noindent {\bf Abstract:}
Kelly's combinatorial lemma is a basic tool in the study of Ulam's reconstruction conjecture. A generalization in terms of a family of $t$-elements subsets of a $v$-element set was given by Pouzet. We consider a version of this generalization modulo a prime $p$.
We give illustrations to graphs and tournaments.\\




\section{Introduction} \label{section def}
Kelly's combinatorial lemma is the assertion that the number $s(F,G)$ of induced subgraphs  of a given graph $G$, isomorphic to $F$, is determined by the deck of $G$, provided that $\vert V(F)\vert < \vert V(G)\vert$, namely $s(F,G) = \frac{1}{\vert V(G)\vert  - \vert V(F)\vert} \sum_{x\in  V(G)} s(F,G_{-x})$ (where $G_{-x}$ is the graph induced by $G$ on $ V(G)\setminus \{x\}$).\\
In terms of a family $\mathcal F$ of $t$-elements subsets of a  $v$-element set, it simply says that
$\vert \mathcal F \vert =  \frac{1}{v-t} \sum_{x\in  V(G)} \vert \mathcal F_{-x} \vert$ where
$ \mathcal F_{-x}:=  \mathcal F \cap [E\setminus \{x\}]^t$. \\
Pouzet  \cite{Pm,Pm2} gave the following extension of this result.
\begin{lemma} (M. Pouzet  \cite{Pm})  \label{lem po}
Let $t$ and $r$ be integers, $V$ be a set of size $v\geq t+r$ elements,
$U$ and $U'$ be sets of subsets $T$ of $t$ elements of $V$. If for every subset
$K$ of $k=t+r$ elements of $V$, the number of elements of $U$ which are contained
in $K$ is equal to the number of elements of $U'$ which are contained in $K$,
then for every finite subsets $T'$ and $K'$ of $V$, such that $T'$ is
contained in $K'$ and  $K'\setminus T'$ has at least $t+r$ elements, the
number of elements of $U$ which contain $T'$ and are contained in $K'$ is equal
to the number of elements of $U'$ which contain $T'$ and are contained in $K'$.
\end{lemma}

In particular if $\vert V \vert \geq 2t+r=t+k$, we have this   particular  version of the combinatorial lemma of Pouzet :

\begin{lemma}  (M. Pouzet \cite{Pm})  \label{particular mp}
Let $v,t$ and $k$ be integers, $V$ be a set of $v$ elements with $t\leq min{(k,v-k)}$,
$U$ and $U'$ be sets of subsets $T$ of $t$ elements of $V$. If for every subset
$K$ of $k$ elements of $V$, the number of elements of $U$ which are contained
in $K$ is equal to the number of elements of $U'$ which are contained in $K$,
then  $U=U'$.
\end{lemma}

We denote by $n(U,K)$ the  number of elements of $U$ which are contained
in $K$, thus Lemma \ref{particular mp} says that if $n(U,K)=n(U',K)$  for every subset
$K$ of $k$ elements of $V$ then $U=U'$.
Here we consider the case where $n(U,K)\equiv n(U',K)$ modulo a prime $p$  for every subset
$K$ of $k$ elements of $V$; our main result, Theorem  \ref{thm js}, is then a version, modulo a prime $p$, of the particular version of the combinatorial lemma  of  Pouzet.\\

Kelly's combinatorial lemma is a basic tool in the study of Ulam's reconstruction conjecture.
Pouzet's combinatorial lemma has been used several times in reconstruction problems (see for example \cite{ ABB, B, BD, BL,  D1, D2}). Pouzet gave a proof of his lemma via a counting argument \cite{Pm2} and latter by using linear algebra (related to incidence matrices) \cite{Pm} (the paper was published earlier).



Let $n,p$ be positive integers, the decomposition of $n=\sum_{i=0}^{n(p)} n_i p^i$ in the basis $p$ is also denoted $[n_0,n_1,\dots ,n_{n(p)}]_p$  where $n_{n(p)}\neq 0$ if and only if $n\neq 0$.

\begin{theorem} \label{thm js}
Let $p$ be a prime number. Let $v,t$ and $k$ be non-negative integers,
$k=[k_0,k_1,\dots , k_{k(p)}]_p$, $t=[t_0,t_1,\dots , t_{t(p)}]_p$. Let
$V$ be a set of $v$ elements with $t\leq min{(k,v-k)}$,
$U$ and $U'$ be sets of subsets $T$ of $t$ elements of $V$.
We assume that for every subset
$K$ of $k$ elements of $V$, the number of elements of $U$ which are contained
in $K$ is equal (mod $p$) to the number of elements of $U'$ which are contained in $K$.\\
1)  If $k_i=t_i$ for all $i<t(p)$ and $k_{t(p)}\geq t_{t(p)}$,  then  $U=U'$.\\
2)  If  $t=t_{t(p)}p^{t(p)}$ and $k=\sum_{i={t(p)}+1}^{k(p)} k_{i}p^i$, we have  $U=U'$, or one of the sets $U,U'$ is the set of all $t$ element-subsets of $V$ and the other is empty, or (whenever $p=2$) for all $t$-element subsets $T$ of $V$,  $T\in U$ if and only if  $T\not\in U'$.
\end{theorem}

Our proof of Theorem \ref{thm js}  is  an application of properties of incidence matrices  due to D.H. Gottlieb \cite{Go}, W. Kantor \cite{KA} and R.M. Wilson \cite{W}, we use Wilson's Theorem  (Theorem \ref{thm Wilson}).  \\
 In  a reconstruction problem  of graphs up to complementation \cite{dlps1}, Wilson's Theorem yielded the following  result:
\begin{theorem} (\cite{dlps1})\label{k=0[4],p=2}
Let  $k$ be an integer,  $2\leq k\leq v-2$, $k\equiv 0$ (mod $4$). Let $G$ and $G'$ be two graphs on the same set $V$ of $v$
vertices (possibly infinite). We assume that $e(G_{\restriction K})$ has the same parity as   $e(G'_{\restriction K})$   for  all $k$-element subsets $K$ of $V$. Then $G'=G$ or $G'=\overline {G}$.
\end{theorem}

Here we look for similar results whenever $e(G_{\restriction K}) \equiv e(G'_{\restriction K})$ modulo a prime  $p$. As an illustration of Theorem \ref{thm js}, we obtain the following result.

\begin{theorem}\label{k=2[4]}
Let  $p$ be a prime number and $k$ be an integer,  $2\leq k\leq v-2$. Let $G$ and $G'$ be two graphs on the same set $V$ of $v$
vertices (possibly infinite). We assume that for all k-element subsets $K$ of $V$, $e(G_{\restriction K}) \equiv e(G'_{\restriction K})$ (mod $p$).\\
1)  If $p\geq3$, $k\not\equiv 0,1 \ (mod \ p)$, then  $G'=G$.\\
2)  If $p\geq3$, $k\equiv 0$ (mod $p$),
then  $G'=G$, or one of the graphs  $G,G'$ is the complete graph and the other is the empty graph.\\
3)   If $p=2$, $k\equiv 2$ (mod $4$), then $G'=G$.
\end{theorem}

We give another  illustrations of Theorem \ref{thm js}, to graphs in section \ref{section graphs}, and to tournaments
 in section \ref{section tournaments}.

\section{Incidence matrices}

We consider the  matrix $W_{t\;k}$ defined as follows :
Let $V$ be a finite set, with $v$ elements. Given non-negative integers $t,k$, let $W_{t\;k}$ be the  ${v \choose t}$  by  ${v \choose k}$ matrix of $0$'s and $1$'s, the rows of which are indexed by the $t$-element subsets
 $T$ of $V$, the columns are indexed by the $k$-element subsets $K$
of $V$, and where the entry $W_{t\;k}(T,K)$ is $1$ if
 $T\subseteq K$ and is $0$ otherwise. The matrix transpose of   $W_{t\;k}$ is  denoted $^tW_{t\;k}$.\\
We say that a matrix $D$ is a {\it {diagonal form}} for a matrix $M$ when $D$ is diagonal and there exist unimodular matrices (square integral matrices which have integral inverses) $E$ and $F$ such that $D=EMF$. We do not require that $M$ and $D$ are square; here "diagonal" just means that the $(i,j)$ entry of $D$ is $0$ if $i\neq j$.
A fundamental result, due to R.M.Wilson \cite{W}, is the following.
\begin{theorem} (R.M. Wilson \cite{W}) \label{thm Wilson+} For $t\leq min{(k,v-k)}$,  $W_{t\; k}$ has as a diagonal form the ${v  \choose t}\times {v  \choose k}$ diagonal matrix with diagonal entries
$$ {k-i  \choose t-i}\ \mbox{with multiplicity}\   {v  \choose i} -   {v  \choose i-1}, \ \ \ i=0,1,\dots ,t.$$
\end{theorem}

Clearly from Theorem \ref{thm Wilson+}, $rank\ W_{t\; k}$  over the field  $\Q$ is ${v  \choose t}$, that is
Theorem \ref{gottlieb-kantor} due to Gottlieb  \cite{Go}. On the other hand, from  Theorem \ref{thm Wilson+}, follows  $rank\ W_{t\; k}$  over the field $\Z/p\Z$, as given by Theorem \ref{thm Wilson}.

\begin{theorem} (R.M. Wilson \cite{W}) \label{thm Wilson} For $t\leq min{(k,v-k)}$, the rank of $W_{t\; k}$ modulo a prime $p$ is
$$ \sum  {v  \choose i}-  {v\choose
i - 1}
$$
where the sum is extended over those indices $i$, $0\leq i\leq t$,  such that $p$
does not divide the binomial coefficient ${k-i \choose
t-i}$.
\end{theorem}

In the statement of the theorem,  ${v \choose -1}$
should be interpreted as zero.\\

A fundamental result, due to D.H. Gottlieb \cite{Go}, and independently W. Kantor \cite {KA}, is this:
\begin{theorem} (D.H. Gottlieb \cite{Go}, W. Kantor \cite {KA}) \label{gottlieb-kantor} For $t\leq min{(k,v-k)}$,  $W_{t\; k}$ has full row rank over the field $\Q$ of rational numbers.
 \end{theorem}

It is clear that  $t\leq min{(k,v-k)}$ implies ${v  \choose t}\leq {v  \choose k}$ then, from Theorem \ref{gottlieb-kantor},  we have the following result :

\begin{corollary}\label{rk-gottlieb-kantor} For $t\leq min{(k,v-k)}$, the  rank of
$W_{t\; k}$ over the field $\Q$ of rational numbers is ${v  \choose t}$  and thus $Ker(^tW_{t\; k})=\{0\}$.
 \end{corollary}

 If $k:=v-t$ then, up to a relabelling, $W_{t\; k}$ is the adjacency matrix $A_{t,v}$
 of the  {\it Kneser graph} $KG(t,v)$ \cite{GoRo}, graph  whose vertices are the $t$-element
subsets of $V$, two subsets forming an edge if they are disjoint.
The eigenvalues of Kneser graphs are computed in \cite{GoRo} (Theorem 9.4.3), and thus an equivalent form of Theorem \ref{gottlieb-kantor}  is:
\begin{theorem}  \label{Ka} $A_{t, v}$ is non-singular for $t\leq \frac{v}{2}$.
\end{theorem}

%

We characterize values of $t$ and $k$ so that $dim \ Ker(^tW_{t\; k})\in \{0,1\}$
 and give a basis of  $Ker(^tW_{t\; k})$, that appears in the following result.

\begin{theorem} \label{thm js2}
Let $p$ be a prime number. Let $v,t$ and $k$ be non-negative integers,
$k=[k_0,k_1,\dots , k_{k(p)}]_p$, $t=[t_0,t_1,\dots , t_{t(p)}]_p$,  $t\leq min{(k,v-k)}$.
We have:\\
1)  $k_j=t_j$ for all $j<t(p)$ and $k_{t(p)}\geq t_{t(p)}$ if and only if   $Ker(^tW_{t\; k})=\{0\}$  \mbox{(mod $p$)}.\\
2)  $t=t_{t(p)}p^{t(p)}$ and $k=\sum_{i={t(p)}+1}^{k(p)}k_ip^i$   if and only if
$\dim Ker (^tW_{t\; k})= 1$\ \mbox{(mod $p$)} and $\{(1,1,\cdots ,1)\}$ is a basis of  $Ker (^tW_{t\; k})$.
\end{theorem}

The proof of  Theorem \ref{thm js2} uses  Lucas's Theorem.
The notation $a\mid b$ (resp. $a\nmid b$) means $a$ divide $b$ (resp. $a$ not divide $b$).

\begin{theorem} (Lucas's Theorem  \cite{Lucas})    \label{lucas}
Let $p$ be a prime number, $t,k$ be   positive integers,
$t\leq k$, $t=[t_0,t_1,\dots ,t_{t(p)}]_p$ and $k=[k_0,k_1,\dots ,k_{k(p)}]_p$. Then
$${k  \choose t} = \prod_{i=0}^{t(p)} {k_i  \choose t_i} \  (mod \ p),\ \mbox{where}  \
{k_i  \choose t_i} =0\ \mbox{if} \ t_i>k_i.$$
\end{theorem}

As a consequence of Theorem \ref{lucas},  we have the following result which is very useful in this paper.

\begin{corollary}\label{cor-lucas}
Let $p$ be a prime number, $t,k$ be   positive integers,
$t\leq k$, $t=[t_0,t_1,\dots ,t_{t(p)}]_p$ and $k=[k_0,k_1,\dots ,k_{k(p)}]_p$. Then
\begin{center}
 $p \vert {k  \choose t}$ if and only if there is $i\in \{0,1,\dots , t(p)\}$ such that $t_i>k_i$.
\end{center}
\end{corollary}

\noindent{\bf{Proof of  Theorem \ref{thm js2}.}} 1) We begin by the direct implication.
We will prove $p \nmid {{k-i}  \choose {t-i}}$ for all $i=[i_0,i_1,\dots , i_{t(p)}]\in \{0,\dots ,t\}$ with
$i_{t(p)}\leq t_{t(p)}$.
Since $k_j=t_j$ for all $j< t(p)$, then  $(t-i)_j=(k-i)_j$ for all $j < t(p)$.
As  $k_{t(p)} \geq  t_{t(p)}\geq i_{t(p)}$ then  $(k-i)_{t(p)} \geq  (t-i)_{t(p)}$,  thus,  by Corollary  \ref{cor-lucas}, $p \nmid {{k-i}  \choose {t-i}}$ for all  $i\in \{0,1, \dots , t\}$.
Now from Theorem \ref{thm Wilson},
$rank \ W_{tk} = \sum_{i=0}^{t} {v  \choose i} - {v  \choose {i-1}}=
{v  \choose t}$. Then  the kernel of $^tW_{t\; k}\ \mbox{(mod $p$)} \ \mbox{is}\ \{0\}$.\\
Now we prove the converse implication.   From Theorem \ref{thm Wilson+},   $Ker(^tW_{t\; k})=\{0\}$ implies  $p \nmid {k-i \choose t-i}$ for all $i\in \{0,1, \dots , t\}$, in particular
$p  \nmid {k \choose t}$.  Then by Corollary   \ref{cor-lucas},
 $k_j\geq t_j$ for all $j \leq t(p)$.
We will prove that  $k_j =  t_j$ for all $j \leq t(p)-1$.
By contradiction, let $s$  be the least integer in $\{0,1,  \dots , t(p)-1\}$, such that $k_s>t_s$.
We have
$(t-(t_s+1)p^s)_s = p-1$, $(k-(t_s+1)p^s)_s = k_s-t_s-1$ and $p-1>k_s-t_s-1$. From Corollary   \ref{cor-lucas},
 $p  \mid {{k-(t_s+1)p^s}  \choose {t-(t_s+1)p^s}}$, that is impossible.\\
2)  Set $n:=t(p)$.  We begin by the direct implication.   Since $0=k_n<t_n$ then, by   Corollary   \ref{cor-lucas},  $p \vert {{k}  \choose {t}}$.
We will prove $p  \nmid {{k-i}  \choose {t-i}}$ for all $i=[i_0,i_1,\dots , i_{n}]\in \{1,2,\dots ,t\}$. \\
Since $k_j=t_j=0$ for all $j<n$, then  $(t-i)_j=(k-i)_j$ for all $j < n$.
From $t_n\geq i_n$, we have  $(t-i)_n\in \{t_n- i_n,t_n- i_n-1\}$.
Note that  $(k-i)_n\in\{p- i_n-1,p- i_n\}$ and $p-i_n-1\geq t_n-i_n$; thus  $(k-i)_n \geq  (t-i)_n$. So for all $j\leq n$,
$(k-i)_j \geq  (t-i)_j$. Then, by Corollary  \ref{cor-lucas},  $p \nmid {{k-i}  \choose {t-i}}$ for all $i\in \{1,2,\dots ,t\}$. Now from Theorem \ref{thm Wilson},
$rank \ W_{tk} = \sum_{i=1}^{t} {v  \choose i} - {v  \choose {i-1}}=
{v  \choose t}-1$, and thus $\dim Ker (^tW_{t\; k})= 1$. Now
  $(1,1,\cdots ,1)W_{t\;  k}=({k  \choose t},{k  \choose t},\cdots ,{k  \choose t})$.\\
Since  $p \mid {k  \choose t}$, then
 $(1,1,\cdots ,1)W_{t\;  k}\equiv0$ (mod $p$).
Then $\{(1,1,\cdots ,1)\}$ is a basis of  the kernel of $^tW_{t\; k}$ (mod $p$).\\
Now we prove the converse implication. Since  $\{(1,1,\cdots ,1)\}$  is a basis of  the kernel of $^tW_{t\; k}$ (mod $p$) and
$(1,1,\cdots ,1)W_{t\;  k}=({k  \choose t},{k  \choose t},\cdots ,{k  \choose t})$, then $p \mid {k  \choose t}$. Since
$dim \ Ker(^tW_{t\; k})=1$, then from Theorem \ref{thm Wilson},
$p  \nmid {k-i \choose t-i}$ for all $i\in \{1,2,\dots ,t\}$.\\
First, let us prove that $t=t_np^n$. Note that $t_n\neq 0$ since $t\neq 0$.
Since $p \vert {k \choose t}$ then, from Corollary   \ref{cor-lucas}, there is an integer $j\in \{0,1,\dots ,n\}$ such that
$t_j > k_j$. Let $A:=\{ j<n\ : \ t_j\neq 0\}$. By contradiction,
assume $A\neq \emptyset$. \\
Case 1. There is $j\in A$ such that $t_j > k_j$. We have $(t-p^n)_j = t_j$, $ (k-p^n)_j=k_j$.
Then from Corollary   \ref{cor-lucas}, we have  $p  \mid {{k-p^n}  \choose {t-p^n}}$, that is impossible.\\
Case 2. For all $j\in A$, $t_j \leq  k_j$. Then  $t_n > k_n$.
We have $(t-p^j)_n = t_n$, $ (k-p^j)_n=k_n$.
Then, from Corollary   \ref{cor-lucas}, we have  $p  \mid {{k-p^j}  \choose {t-p^j}}$, that is impossible.\\
From the above  two cases, we deduce $t=t_np^n$.\\
Secondly, since  $p \vert {{k}  \choose {t}}$, then by  Corollary   \ref{cor-lucas}, $t_n>k_n$.
Let us show that $k_n=0$.
By contradiction, if $k_n\neq 0$ then
$(t-p^n)_n=t_{n}-1>  k_n-1=(k-p^n)_n$. From Corollary   \ref{cor-lucas},
$p  \mid {{k-p^n}  \choose {t-p^n}}$, that is impossible.
Let $s\in \{0,1,\dots ,n-1\}$, let us show that $k_s=0$.
By contradiction, if $k_s\neq 0$ then, $(t-p^s)_s =p-1$, $(k-p^s)_s = k_s-1$, thus    $(t-p^s)_s > (k-p^s)_s$
so,  from Corollary   \ref{cor-lucas},
 $p  \mid {{k-p^s}  \choose {t-p^s}}$, that is impossible.\endproof

%

\section{Proof of  Theorem \ref{thm js}.}

Let $T_1,T_2, \cdots  ,T_{{v \choose  t}}$  be an enumeration of the   $t$-element subsets of $V$,  let $K_1,K_2, \cdots ,K_{{v \choose  k}}$ be an enumeration of the   $k$-element
subsets of $V$ and  $W_{t\; k}$ be the matrix of the $t$-element subsets versus the $k$-element subsets.

Let $w_U$ be the row matrix $(u_1,u_2, \cdots , u_{v \choose  t})$ where $u_i=1$ if   $T_i\in U$,  $0$ otherwise.   We have
   $$w_UW_{t\;  k}=(\vert \{T_i\in U : T_i \subseteq K_1\}\vert , \cdots ,\vert \{ T_i\in U : T_i \subseteq K_{{v \choose  k}} \}\vert).$$
$$w_{U'}W_{t\;  k}=(\vert \{T_i\in U' : T_i \subseteq K_1\}\vert , \cdots ,\vert \{ T_i\in U' : T_i \subseteq K_{{v \choose  k}} \}\vert).$$

Since for all $j\in \{1,\dots ,{v \choose  k}\}$, the number of elements of $U$ which are contained in $K_j$  is equal (mod $p$) to the number of elements of $U'$ which are contained in $K_j$, then    $(w_U-w_{U'})W_{t\; k}=0$ \ \mbox{(mod $p$)}, so
$w_U-w_{U'}\in Ker (^tW_{t\; k})$.\\
1) Assume  $k_i=t_i$ for all $i<t(p)$ and $k_{t(p)}\geq t_{t(p)}$. From 1) of Theorem \ref{thm js2}, $w_U-w_{U'}=0$, that gives $U=U'$.\\
2) Assume  $t=t_{t(p)}p^{t(p)}$ and $k=\sum_{i={t(p)}+1}^{k(p)}k_ip^i$. From 2) of Theorem  \ref{thm js2}, there is an integer $\lambda \in [0,p-1]$ such that $w_U-w_{U'}=\lambda (1,1,\cdots ,1)$. It is clear that  $\lambda \in \{0,1,-1\}$.
 If  $\lambda =0$ then $U=U'$. If $\lambda =1$ and $p\geq 3$ then
$U=\{T_1,T_2, \cdots  ,T_{{v \choose  t}}\}$, $U'=\emptyset$.  If $\lambda =1$ and $p= 2$ then
$U=\{T_1,T_2, \cdots  ,T_{{v \choose  t}}\}$, $U'=\emptyset$, or $T\in U$  if and only if  $T\not\in U'$.
 If $\lambda =-1$ and $p\geq 3$ then $U=\emptyset$,
$U'=\{T_1,T_2, \cdots  ,T_{{v \choose  t}}\}$.
 If $\lambda =-1$ and $p= 2$ then
$U'=\{T_1,T_2, \cdots  ,T_{{v \choose  t}}\}$, $U=\emptyset$, or $T\in U$  if and only if   $T\not\in U'$.
\endproof

\section{Illustrations to graphs} \label{section graphs}
Our notations and terminology follow \cite {Bo}.
A \textit{digraph}  $G = (V, E)$ or $G=(V(G),E(G))$, is formed
 by a finite set $V$ of vertices and  a set $E$ of pairs of distinct vertices, called {\it arcs} of $G$. The  {\it order} (or {\it cardinal}) of $G$ is the  number of its vertices.  If $K$ is a subset of $V$, the {\it  restriction} of $G$ to $K$, also called the {\it induced subdigraph} of $G$ on $K$ is the digraph $G_{\restriction K}:= (K, K^2\cap  E)$. If $K=V\setminus \{x\}$, we denote this  digraph by $G_{-x}$. Let $G = (V, E)$ and $G' = (V', E')$ be two digraphs. A one-to-one correspondence $f$ from $V$ onto $V'$ is an {\it isomorphism from} $G$ {\it onto} $G'$ provided that for $x, y \in V$, $(x, y) \in E$ if and only if $(f(x), f(y)) \in E'$. The digraphs $G$ and $G'$ are then said to be {\it isomorphic}, which is denoted by  $G \simeq G^{\prime}$. A subset $I$ of $V$ is an
\textit{interval} \cite{FR,ST} (or a {\it clan}  \cite{ER}, or an {\it homogenous subset} \cite{TG}) of  $G$ provided that for all $a, b\in
I$ and  $x \in V\setminus I$, $(a, x) \in E(G)$  if and only if   $(b, x)\in  E(G)$,  and the same for $(x,a)$ and $(x,b)$. For example $\emptyset$, $\{x\}$
where $x \in V$, and $V$ are intervals of $G$, called {\it trivial intervals}. A digraph is then said to be {\it indecomposable} \cite{ST} (or {\it primitive }\cite{ER}) if all its intervals are trivial, otherwise it is said to be {\it decomposable}.\\
We say that $G$ is a {\it graph} (resp. {\it tournament}) when for every distinct vertices $x,y$ of $V$, $(x, y) \in E$ if and only if $(y, x)\in E$ (resp $(x, y) \in E$ if and only if $(y, x)\not\in E$); we say that $\{x,y\}$ is an {\it edge} of the graph $G$ if $(x,y)\in E$, thus $E$ is identified with a subset of $[V]^2$, the set of pairs $\{x,y\}$ of distinct elements of $V$.\\
Let $G= (V, E)$ be a graph,
the {\it complement} of  $G$ is the graph  $\overline G:= (V, [V]^2\setminus E)$. We denote by $e(G):=\vert E(G)\vert $ the number of edges of  $G$.
The {\it degree} of a vertex $x$
of $G$, denoted $d_G(x)$, is  the number of  edges which contain $x$.
A $3$-element subset  $T$ of $V$ such that all pairs belong to $E(G)$ is a  {\it triangle} of $G$.
Let $T(G)$ be the set
of {\it triangles} of $G$ and let $t(G):=\mid T(G)\mid$.
A $3$-element subset of $V$ which is a triangle of $G$ or of $\overline G$ is a  $3$-{\it homogeneous} subset of $G$.
We set
$H^{(3)}(G):=T(G)\cup T(\overline G)$,  the set of $3$-homogeneous subsets of $G$, and $h^{(3)}(G):=\mid H^{(3)}(G)\mid$.\\

\noindent{\bf Another proof of Theorem \ref{k=0[4],p=2} using Theorem \ref{thm js}.} Here 
$p=2$, $t=2=[0,1]_p$ and $k=[0,0,k_2,\dots]_p$. From 2) of  Theorem \ref{thm js}, $U=U'$, or one of the sets $U,U'$ is the set of all $2$ element-subsets of $V$ and the other is empty, or for all $2$-element subsets $T$ of $V$,  $T\in U$  if and only if  $T\not\in U'$. Thus $G'=G$ or $G'=\overline{G}$.\endproof
\\

\noindent{\bf Proof of Theorem \ref{k=2[4]}.}
We set $U:=E(G)$,  $U':=E(G')$.
For all $K\subseteq V$ with $\vert K\vert = k$, we have:
$\{\{x,y\}\subseteq K : \{x,y\} \in U\}= E(G_{\restriction K})$ and
$\{\{x,y\}\subseteq K : \{x,y\} \in U'\}= E(G'_{\restriction K})$.
Since $e(G_{\restriction K}) \equiv e(G'_{\restriction K})$ (mod $p$), then
$\vert\{\{x,y\}\subseteq K : \{x,y\} \in U\}\vert   \equiv \vert\{\{x,y\}\subseteq K : \{x,y\} \in U'\}\vert$ (mod $p$).\\
1) $p\geq3$, $t=2=[2]_p$ and $k_0\geq 2$. From 1) of  Theorem \ref{thm js},  $U=U'$, thus $G=G'$.\\
2) $p\geq3$, $t=2=[2]_p$ and $k_0=0$. From 2) of  Theorem \ref{thm js},  we have $U=U'$ or
one of $U,U'$ is the set of all $2$-elements subsets of $V$ and the other is empty.
Then $G=G'$ or  one of the graphs  $G,G'$ is the complete graph and the other is the empty graph.\\
3)  $p=2$, $t=2=[0,1]_p$ and $k=[0,1,k_2,\dots]_p$. From 1) of  Theorem \ref{thm js}, we have $U=U'$, thus $G=G'$.
\endproof
\\

The following result  concerns  graphs $G$ and $G'$ such that  $h^{(3)}(G_{\restriction K})\equiv h^{(3)}(G'_{\restriction K})$ modulo a prime $p$,
for all $k$-element subsets $K$ of $V$.

\begin{theorem}  \label{Ka+lem+TR} Let $G$ and $G'$ be two graphs on the
same set $V$ of $v$ vertices.  Let $p$ be a prime number  and $k$ be an integer, $3\leq k\leq v-3$.\\
1)   If $h^{(3)}(G_{\restriction K})=h^{(3)}(G'_{\restriction K})$
for all $k$-element subsets $K$ of $V$ then $G$ and $G'$ have the same $3$-element homogeneous sets.\\
2)  Assume $p\geq 5$.  If $k\not\equiv 1,2 \ (mod \ p)$ and  $h^{(3)}(G_{\restriction K})\equiv h^{(3)}(G'_{\restriction K})$ (mod $p$)
for all $k$-element subsets $K$ of $V$, then $G$ and $G'$ have the same $3$-element homogeneous sets.\\
3)  If ($p=2$ and $k\equiv 3 \ (mod \ 4)$) or   ($p=3$ and $3\mid k$), and  $h^{(3)}(G_{\restriction K})\equiv h^{(3)}(G'_{\restriction K})$   (mod $p$)
for all $k$-element subsets $K$ of $V$, then $G$ and $G'$ have the same $3$-element homogeneous sets.
\end{theorem}

\Proof  $H^{(3)}(G)=\{\{a,b,c\}: G_{\restriction \{a,b,c\}} \ \mbox{is a $3$-element homogeneous set}\}$.
We set $U:=H^{(3)}(G)$ and $U':=H^{(3)}(G')$.
For all $K\subseteq V$ with $\vert K\vert = k$, we have:
$\{T\subseteq K : T \in U\}=H^{(3)}_{G_{\restriction K}}$ and
$\{T\subseteq K : T \in U'\}=H^{(3)}_{G'_{\restriction K}}$. Set $t:=\mid T\mid =3$.\\
1) Since $h^{(3)}(G_{\restriction K})=h^{(3)}(G'_{\restriction K})$
for all $k$-element subsets $K$ of $V$ then
$\vert\{T\subseteq K : T \in U\}\vert   = \vert\{T\subseteq K : T \in U'\}\vert$. From Lemma \ref{particular mp} it follows that $U=U'$, then
 $G$ and $G'$ have the same $3$-element homogeneous sets. \\
2) Since $h^{(3)}(G_{\restriction K})\equiv h^{(3)}(G'_{\restriction K})$ (mod $p$)
for all $k$-element subsets $K$ of $V$ then
$\vert\{T\subseteq K : T \in U\}\vert   \equiv \vert\{T\subseteq K : T \in U'\}\vert$ (mod $p$). \\
Case 1. $p\geq5$, $t=3=[3]_p$, $k=[k_0,\dots]_p$ and $t_0=3\leq k_0$. From 1) of Theorem  \ref{thm js} we have $U=U'$, thus  $G$ and $G'$ have the same $3$-element homogeneous sets.\\
Case 2. $p\geq5$, $t=3=[3]_p$, $k=[0,k_1,\dots]_p$.
By Ramsey's Theorem \cite {Ra},  every graph with at least $6$ vertices contains
a  $3$-element homogeneous set. Then $U$ and $U'$ are nonempty, so  from  2) of Theorem  \ref{thm js}, $U=U'$, thus
$G$ and $G'$ have the same $3$-element homogeneous sets.\\
3) Since $h^{(3)}(G_{\restriction K})\equiv h^{(3)}(G'_{\restriction K})$ (mod $p$)
for all $k$-element subsets $K$ of $V$ then
$\vert\{T\subseteq K : T \in U\}\vert   \equiv \vert\{T\subseteq K : T \in U'\}\vert$ (mod $p$). \\
Case 1. $p=2$, $t=3=[1,1]_p$ and $k\equiv 3 \ (mod \ 4)$. In this case, $k=[1,1,k_2, \dots]_p$, then
 from 1) of Theorem  \ref{thm js} we have $U=U'$, thus  $G$ and $G'$ have the same $3$-element homogeneous sets.\\
Case 2. $p=3$, $t=3=[0,1]_p$ and $k=[0,k_1,\dots,k_{k(p)}]_p$.\\
Case 2.1. $k_1\in\{1,2\}$, then from 1) of Theorem  \ref{thm js} we have $U=U'$, thus  $G$ and $G'$ have the same $3$-element homogeneous sets.\\
Case 2.2. $k_1=0$.
By Ramsey's Theorem \cite {Ra},  every graph with at least $6$ vertices contains
a  $3$-element homogeneous set. Then $U$ and $U'$ are nonempty, so  from  2) of Theorem  \ref{thm js}, $U=U'$, thus
$G$ and $G'$ have the same $3$-element homogeneous sets.
\endproof
\\

Let $G=(V,E)$ be a graph.
 From \cite{ST}, every  indecomposable graph of size $4$ is isomorphic to  $P_4=\left(\{0,1,2,3\},\{\{0,1\},\{1,2\},\{2,3\}\}\right)$. Let ${\mathcal P}^{(4)}(G)$ be the set of indecomposable induced subgraphs of $G$ of size $4$,
we set  $p^{(4)}(G):=\vert{\mathcal P}^{(4)}(G)\vert$. The following result  concerns  graphs $G$ and $G'$ such that $p^{(4)}(G_{\restriction K})\equiv p^{(4)}(G'_{\restriction K})$ modulo a prime $p$,
for all $k$-element subsets $K$ of $V$.

\begin{theorem}  \label{Ka+lem+P4} Let $G$ and $G'$ be two graphs on the
same set $V$ of $v$ vertices. Let $p$ be a prime number and $k$  be an integer,  $4\leq k\leq v-4$.\\
1)   If $p^{(4)}(G_{\restriction K})=p^{(4)}(G'_{\restriction K})$
for all $k$-element subsets $K$ of $V$ then $G$ and $G'$ have the same indecomposable sets of size $4$.\\
2)  Assume  $p^{(4)}(G_{\restriction K})\equiv p^{(4)}(G'_{\restriction K})$  (mod $p$)
for all $k$-element subsets $K$ of $V$.\\
a)  If $p\geq 5$ and  $k\not\equiv 1,2,3 \ (mod \ p)$, then  $G$ and $G'$ have the same indecomposable sets of size $4$.\\
b)  If ($p=2$, $4\mid k$ and $8\nmid k$) or ($p=3$, $3\mid k-1$ and $9\nmid k-1$),  then  $G$ and $G'$ have the same indecomposable sets of size $4$.\\
c)  If  $p=2$ and $8\mid k$, then  $G$ and $G'$ have the same indecomposable sets of size $4$, or for all $4$-element subsets  $T$ of $V$, $G_{\restriction T}$ is indecomposable if and only if $G'_{\restriction T}$ is decomposable.
\end{theorem}

\Proof
Let $U:=\{T\subseteq V : \vert T\vert = 4, \ G_{\restriction T}\simeq P_4 \}={\mathcal P}^{(4)}(G) $,
$U':=\{T\subseteq V : \vert T\vert = 4, \ G'_{\restriction T}\simeq P_4 \}={\mathcal P}^{(4)}(G') $.
For all $K\subseteq V$, we have
$\{T\subseteq K : T\in U\}= {\mathcal P}_4( G_{\restriction K})$ and
$\{T\subseteq K : T\in U'\}= {\mathcal P}_4( G'_{\restriction K})$. Set $t:= \vert T \vert =4$.\\
1)  Since $p^{(4)}(G_{\restriction K})=p^{(4)}(G'_{\restriction K})$ then
$\vert\{T\subseteq K : T\in U\}\vert=\vert\{T\subseteq K : T\in U'\}\vert$.
From Lemma \ref{particular mp},  $U=U'$, then $G$ and $G'$ have the same indecomposable sets of size $4$.\\
2) We have $p^{(4)}(G_{\restriction K})\equiv p^{(4)}(G'_{\restriction K})$ (mod $p$)
for all $k$-element subsets $K$ of $V$, then $\vert\{T\subseteq K : T\in U\}\vert \equiv \vert\{T\subseteq K : T\in U'\}\vert$ (mod $p$).\\
a) Case 1. $p\geq5$, $t=4=[4]_p$, $k=[k_0,\dots]_p$ and $t_0=4\leq k_0$. From 1) of Theorem  \ref{thm js} we have $U=U'$, thus $G$ and $G'$ have the same indecomposable sets of size $4$.\\
Case 2. $p\geq5$, $t=4=[4]_p$, $k=[0,k_1,\dots]_p$. Since in every graph of order $5$, there is a restriction of size $4$ not isomorphic to $P_4$ then, from  2) of Theorem  \ref{thm js}, $U=U'$,  thus
$G$ and $G'$ have the same indecomposable sets of size $4$.\\
b)
Case 1. $p=2$, $t=4=[0,0,1]_p$ and $k=[0,0,1,k_3,\dots,k_{k(p)}]_p$.
From 1) of Theorem  \ref{thm js},  we have $U=U'$, thus $G$ and $G'$ have the same indecomposable sets of size $4$.\\
Case 2. $p=3$, $t=4=[1,1]_p$, $k=[1,k_1,\dots,k_{k(p)}]_p$ and $t_1=1\leq k_1$. From 1) of Theorem  \ref{thm js},   $U=U'$, thus $G$ and $G'$ have the same indecomposable sets of size $4$.\\
c) We have $p=2$, $t=4=[0,0,1]_p$, $k=[0,0,0,k_3,\dots,k_{k(p)}]_p$. Since in every graph of order $5$, there is a restriction of size $4$  not isomorphic to $P_4$, then from  2) of Theorem  \ref{thm js}, $U=U'$, or for all $4$-element subsets $T$ of $V$,  $T\in U$  if and only if  $T\not\in U'$.  Thus
$G$ and $G'$ have the same indecomposable sets of size $4$, or for all $4$-element subsets $T$ of $V$, $G_{\restriction T}$ is indecomposable if and only if $G'_{\restriction T}$ is decomposable.
\endproof
\\


 In  a reconstruction problem  of graphs up to complementation \cite{dlps1}, Wilson's Theorem yielded the following  result:


 \begin{theorem} (\cite{dlps1})\label{k=1[4]}
Let $G$ and $G'$ be two graphs on the same set $V$ of $v$
vertices (possibly infinite). Let  $k$ be an integer, $5\leq k\leq v-2$,
$k\equiv 1$  (mod $4$). Then the following properties are equivalent:\\
(i)  $e(G_{\restriction K})$ has the same parity as   $e(G'_{\restriction K})$   for  all $k$-element subsets $K$ of $V$;   and   $G_{\restriction K}$, $G'_{\restriction K}$ have the same $3$-homogeneous subsets;\\
(ii) $G'= G$ or $G'= \overline G$.
\end{theorem}


Here, we just want to point out that we can obtain a similar result for $k\equiv 3$ (mod $4$), namely Theorem \ref{k=3[4]},  using the same proof as that of Theorem   \ref{k=1[4]}. 

The {\it boolean sum}  $G\dot{+} G'$
of two graphs $G=(V,E)$ and $G'=(V,E')$  is the graph $U$ on $V$ whose edges are pairs $e$ of
vertices such that $e\in E$ if and only if $e\notin E'$.

 \begin{theorem}\label{k=3[4]}
Let $G$ and $G'$ be two graphs on the same set $V$ of $v$
vertices (possibly infinite). Let  $k$ be an integer, $3\leq k\leq v-2$,
$k\equiv 3$ (mod $4$). Then the following properties are equivalent:\\
(i) $e(G_{\restriction K})$ has the same parity as   $e(G'_{\restriction K})$   for  all $k$-element subsets $K$ of $V$;  and $G_{\restriction K}$, $G'_{\restriction K}$ have the same $3$-homogeneous subsets;\\
(ii) $G'= G$.
\end{theorem}

\Proof
It is exactly the same as that of Theorem   \ref{k=1[4]} (see (\cite{dlps1}). 
 The implication $(ii)\Rightarrow (i)$ is trivial. We prove  $(i)\Rightarrow (ii)$. 
We  suppose $V$ finite, we set $U:= G\dot{+} G'$, let $T_1,T_2, \cdots , T_{{v \choose  2}}$  be an enumeration of the   $2$-element subsets of $V$,  let $K_1,K_2,\cdots ,K_{{v \choose  k}}$ be an enumeration of the   $k$-element
subsets of $V$.
Let $w_U$ be the row matrix $(u_1,u_2,\cdots , u_{v \choose  2})$ where $u_i=1$ if   $T_i$ is an edge of $U$,  $0$ otherwise.\\     We have
   $w_UW_{2\;  k}=(e(U_{\restriction K_1}),e(U_{\restriction K_2}),\cdots ,  e(U_{\restriction K_{v \choose  k}}))$.
 From the facts that $e(G_{\restriction K})$ has the same parity as  $e(G'_{\restriction K})$   and  $e(U_{\restriction K})=e(G_{\restriction K})+e(G'_{\restriction K})-2e(G_{\restriction K}\cap G'_{\restriction K})$ for  all $k$-element subsets $K$, $w_U$ belongs to the kernel of  $^tW_{2\; k}$ over the $2$-element field. According to Theorem \ref{thm Wilson},  the rank of  $W_{2k}$ (mod $2$) is
${v \choose  2} -v+1$. Hence  $\dim Ker(^tW_{2\; k})=v-1$.

We give a similar claim as Claim 2.8 of \cite{dlps1}, the proof is identical.
\begin{claim} \label{bipartite}Let $k$ be an integer such that $3\leq k\leq v-2$,
$k\equiv 3$  (mod $4$),  then the kernel of $^t W_{2\; k}$ consists of complete bipartite graphs (including the empty graph).
\end{claim}
\Proof  Let us recall that a {\it star-graph} of $v$ vertices consists of a vertex linked to all other vertices, those $v-1$ vertices forming an independent set.
First we prove that each star-graph $S$ belongs to $\K$, the kernel of $^t W_{2\; k}$. Let $w_S$ be the row matrix $(s_1,s_2,\cdots , s_{v \choose  2})$ where $s_i=1$ if   $T_i$ is an edge of $S$,  $0$ otherwise.     We have
   $w_SW_{2\;  k}=(e(S_{\restriction K_1}),e(S_{\restriction K_2}),\cdots ,  e(S_{\restriction K_{v \choose  k}}))$. For all $i\in \{1,\dots ,{v \choose  k}\}$, $e(S_{\restriction K_i})=k-1$ if $1\in K_i$, $0$ otherwise. Since $k$ is odd, each star-graph $S$ belongs to $\K$.
   The vector space (over the $2$-element field) generated by the star-graphs on $V$ consists of all complete bipartite graphs; since $v\geq3$, these are distinct from the complete graph (but include  the empty graph). Moreover, its dimension is $v-1$ (a basis being made of star-graphs). Since $\dim Ker(^tW_{2\; k})=v-1$, then $\K$ consists of complete bipartite graphs as claimed.\endproof

A {\it claw} is a star-graph on  four  vertices, that is a graph made of a vertex joined to three other vertices, with no edges between these  three vertices.  A graph is {\it claw-free} if no induced subgraph is a claw.

\begin{claim} \label{clawfree} (\cite{dlps1}) Let $G$ and $G'$ be two graphs on the same set and having the same $3$-homogeneous subsets, then the boolean sum $U: =G\dot {+} G'$ is   claw-free.
\end{claim}

From Claim \ref{bipartite}, $U$  is a complete bipartite graph and,
 from Claim \ref{clawfree},
$U$ is claw-free. Since $v\geq 5$, it follows that $U$ is  the empty graph. Hence $G'=G$ as claimed.
\endproof

\section{Illustrations to tournaments} \label{section tournaments}
Let  $T=(V,E)$ be a tournament.
For two distinct vertices $x$ and $y$ of  $T$, $x\longrightarrow_Ty$ (or simply $x\longrightarrow y$) means that $(x, y)\in E$ and $(y, x)\not\in E$. For $A\subseteq V$ and $y\in V$,  $A\longrightarrow y$ means $x\longrightarrow y$ for all $x\in A$.
The {\it degree} of a vertex $x$
of $T$ is  $d_T(x):=\vert\{   y\in V:x\longrightarrow y\}\vert$. We denote by $T^*$ {\it the dual} of $T$ that is  $T^*=(V,E^*)$ with $(x,y)\in E^*$  if and only if   $(y,x)\in E$.
A {\it transitive} tournament or  a {\it total order } or {\it$k$-chain} (denoted $O_k$) is a tournament of cardinality $k$, such that for $x, y, z \in V$, if $x\longrightarrow y$ and $y\longrightarrow z$, then $x\longrightarrow z$. If $x$ and $y$ are two distinct vertices of a total order, the notation $x < y$ means that $x\longrightarrow y$.
The tournament $C_3 :=\{\{0,1,2\}, \{(0,1),(1,2),(2,0)\}\}$ (resp. $C_4:=(\{0,1,2,3\},\{(0,3), (0,1), (3,1), (1,2), (2,0), (2,3)\})$) is a $3$-{\it cycle} (resp. $4$-{\it cycle}).
A {\it  diamond}  is a
tournament on $4$ vertices admitting only one interval of cardinality $3$ which is a $3$-cycle.
Up to isomorphism, there  are exactly two diamonds $\delta^{+}$ and $\delta^{-}=(\delta^{+})^{*}$, where $\delta^+$ is the tournament defined on $\{0, 1, 2, 3\}$ by $\delta^+_{\restriction\{0, 1, 2\}} = C_3$ and $\{0, 1, 2\}\rightarrow 3 $. A tournament isomorphic to $\delta^+$ (resp. isomorphic to $\delta^-$) is said to be a {\it positive diamond} (resp. {\it negative diamond}).
 The {\it boolean sum}  $U:=T\dot{+} T'$
of two tournaments $T=(V,E)$ and $T'=(V,E')$,   is the graph $U$ on $V$ whose edges are pairs $\{x,y\}$ of
vertices such that $(x,y)\in E$ if and only if $(x,y)\notin E'$.\\


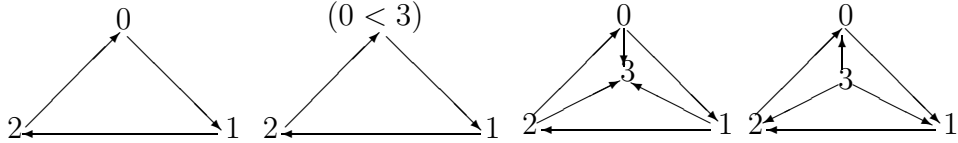
\begin{figure}[H]\label{diamond}
\begin{center}
\setlength{\unitlength}{0.5cm}
\begin{picture}(5,5)

\put(-7.2,5.4){$0$}
\put(-10.1,2.5){$2$} \put(-4.3,2.5){$1$}

\put(-4.5,2.6){\vector(-1,0){5.2}}

\put(-9.6,2.8){\vector(1,1){2.5}}

\put(-6.9,5.2){\vector(1,-1){2.5}}


\put(-1.6,5.5){$(0<3)$}
\put(-3.3,2.5){$2$} \put(2.6,2.5){$1$}

\put(2.5,2.6){\vector(-1,0){5.3}}

\put(-2.7,2.8){\vector(1,1){2.5}}

\put(0,5.2){\vector(1,-1){2.5}}

\put(3.6,2.6){2} \put(8.8,2.6){1}
\put(6.1,5.5){0} \put(6.2,4){3}
\put(8.6,2.9){\vector(-2,1){2.1}}
\put(6.4,5.35){\vector(1,-1){2.4}}
\put(3.9,3.1){\vector(1,1){2.3}}
\put(4,2.9){\vector(2,1){2.2}}
\put(8.6,2.7){\vector(-1,0){4.5}}
\put(6.3,5.4){\vector(0,-1){1}}

\put(9.6,2.6){2} \put(14.8,2.6){1}
\put(12,5.5){0} \put(12.0,3.8){3}
\put(12.4,3.9){\vector(2,-1){2.1}}
\put(12.2,5.4){\vector(1,-1){2.5}}
\put(9.8,3.1){\vector(1,1){2.3}}
\put(12,3.9){\vector(-2,-1){2.0}}
\put(14.6,2.7){\vector(-1,0){4.5}}
\put(12.1,4.3){\vector(0,1){0.9}}
\end{picture}
\vspace{- 1.2cm}
\caption{Cycle $C_3$, Cycle $C_4$, Positive Diamond, Negative Diamond.}
\end{center}
\end{figure}
\vspace{-0.8cm}

\begin{theorem} \label{tournois}Let $T=(V,E)$  and   $T'=(V,E')$ be two tournaments. Let $p$ be a prime number and $k$ be an integer,  $2\leq k\leq v-2$. Let $G:=T \dot{+} T'$. We assume that for all $k$-element subsets  $K$ of $V$, $e(G_{\restriction K})\equiv 0$ (mod $p$).\\
1)  If $p\geq3$, $k  \not\equiv 0,1$ (mod $p$), then   $T'=T$.\\
2)   If $p\geq 3$, $k\equiv 0$ (mod $p$), then $T'=T$ or $T'=T^*$.\\
3)   If $p=2$, $k\equiv 2$ (mod $4$), then $T'=T$.\\
4)   If $p=2$, $k\equiv 0$ (mod $4$), then $T'=T$ or $T'=T^*$.
\end{theorem}

\Proof
We set $G':=$ The empty graph. Then  $e(G_{\restriction K})\equiv e(G'_{\restriction K})$
(mod $p$).\\
1) From 1) of Theorem \ref{k=2[4]},  $G$ is the empty graph,  then  $T'=T$.\\
2) From 2) of Theorem \ref{k=2[4]},  $G$ is empty or the complete graph, then $T'=T$ or $T'=T^*$.\\
3)  From 3) of Theorem \ref{k=2[4]},  $G$ is the empty graph,  then  $T'=T$.\\
4) From Theorem \ref{k=0[4],p=2},  $G$ is the empty graph or the complete graph, then $T'=T$ or $T'=T^*$.
\endproof
\\

Let $T$ be a tournament, we set $C^{(3)}(T):=\{\{a,b,c\} : T_{\restriction \{a,b,c\}} \ \mbox{is a $3$-cycle} \}$, and $c^{(3)}(T):=\mid C^{(3)}(T)\mid$. Let $T=(V,E)$ and $T'=(V,E')$ be two tournaments, let $k$ be a non-negative integer, $T$ and $T'$ are $k$-{\it hypomorphic}  \cite{Bou-Lop,Lr} (resp. $k$-{\it hypomorphic} up to duality)
 if for every $k$-element subset $K$ of $V$, the induced subtournaments $T'_{\restriction K}$ and $T_{\restriction K}$  are isomorphic (resp. $T'_{\restriction K}$ is isomorphic to
$T_{\restriction K}$ or to $T^*_{\restriction K}$). We say that  $T$ and $T'$ are ($\leq k$)-{\it hypomorphic}
 if $T$ and $T'$ are $h$-hypomorphic for every $h\leq k$. Similarly, we say  that  $T$ and $T'$ are  $(\leq k)$-{\it hypomorphic up to duality} if  $T$ and $T'$ are $h$-hypomorphic up to duality for every $h\leq k$.

\begin{theorem}   Let $T$ and $T'$ be two tournaments on the
same set $V$ of $v$ vertices. Let $p$ be a prime number and $k$ be an integer, $3\leq k\leq v-3$.\\
1)  If $c^{(3)}(T_{\restriction K})=c^{(3)}(T'_{\restriction K})$
for all $k$-element subsets $K$ of $V$ then $T$ and $T'$ are $(\leq 3)$-hypomorphic.\\
2) Assume $p\geq5$. If $k\not\equiv 1,2$ (mod $p$), and  $c^{(3)}(T_{\restriction K})\equiv c^{(3)}(T'_{\restriction K})$ (mod $p$)
for all $k$-element subsets $K$ of $V$, then $T$ and $T'$ are $(\leq 3)$-hypomorphic.\\
3)  If ($p=2$ and $k\equiv3$ (mod $4$)) or   ($p=3$ and $3\mid k$), and  $c^{(3)}(G_{\restriction K})\equiv c^{(3)}(G'_{\restriction K})$  (mod $p$)
for all $k$-element subsets $K$ of $V$, then $T$ and $T'$  are $(\leq 3)$-hypomorphic.
\end{theorem}

\Proof
Since every tournament, of cardinality  $\geq 4$, has at least
a restriction of cardinality $3$ which is not a $3$-cycle, then the proof is similar to that  of Theorem \ref{Ka+lem+TR}.\endproof

Let $T$ be a tournament, we set $D^+_4(T):=\{\{a,b,c,d\} : T_{\restriction \{a,b,c,d\}} \simeq \delta^+ \}$,  $D^-_4(T):=\{\{a,b,c,d\} : T_{\restriction \{a,b,c,d\}} \simeq \delta^- \}$, $d^+_4(T):=\mid D^+_4(T)\mid$ and $d^-_4(T):=\mid D^-_4(T)\mid$.

It is well-known that every subtournament of order $4$ of a tournament is either a diamond, a $4$-chain, or a $4$-cycle subtournament. We have $c^{(3)}(O_4)=0$, $c^{(3)}(\delta^+)=c^{(3)}(\delta^-)=1$, $c^{(3)}(C_4)=2$ and $C_4\simeq C_4^*$.
\begin{theorem}\label{tournament}
Let $T$ and $T'$ be two  $(\leq 3)$-hypomorphic tournaments on the
same set $V$ of $v$ vertices. Let $p$ be a prime number and $k$  be an integer,  $4\leq k\leq v-4$.\\
1)  If $d^+_4(T_{\restriction K})=d^+_4(T'_{\restriction K})$
for all $k$-element subsets $K$ of $V$ then $T'$ and $T$ are $(\leq 5)$-hypomorphic.\\
2)  Assume  $d^+_4(T_{\restriction K})\equiv d^+_4(T'_{\restriction K})$ (mod $p$) for all $k$-element subsets $K$ of $V$.\\
 a)
  If $p\geq5$ and,  $k  \not\equiv 1,2,3$ (mod $p$), then $T'$ and $T$ are $(\leq 5)$-hypomorphic.\\
b)   If ($p=3$, $3\mid k-1$ and $9\nmid k-1$) or ($p=2$, $4\mid k$ and $8\nmid k$), then $T'$ and $T$ are $(\leq 5)$-hypomorphic.\\
c)  If $p=2$ and $8\mid k$, then $T'$ and $T$ are $(\leq 5)$-hypomorphic or for all $4$-elements subset $S$ of V, $T_{\restriction S}$ is isomorphic to $\delta^+$ if and only if  $T'_{\restriction S}$ is isomorphic to $\delta^-$.
\end{theorem}

\Proof  To prove that $T'$ and $T$ are ($\leq 5$)-hypomorphic,  the following  lemma shows that it is sufficient to prove that  $T'$ and $T$ are ($\leq 4$)-hypomorphic.
\begin{lemma} \label{hypomorphe} \cite{B}
Let $T$ and $T'$ be two $(\leq 4)$-hypomorphic tournaments
on at least $5$ vertices. Then, $T$ and $T'$ are
$(\leq 5)$-hypomorphic.
\end{lemma}

Now,  let $U^+:=\{S\subseteq V, \ T_{\restriction S} \simeq \delta^+ \}=D^+_4(T) $,
$U'^+:=D^+_4(T') $, $U^-:=D^-_4(T)$ and
$U'^-:=D^-_4(T') $.
\begin{claim}\label{3hyp4hyp}
If $T$ and $T'$ are $(\leq 3)$-hypomorphic and $U^+=U'^+$,  then $U^-=U'^-$; $T$ and $T'$ are $(\leq 5)$-hypomorphic.
\end{claim}
\Proof
Let $S\in U^-$, $T_{\restriction S}\simeq \delta^-$. Since $T$ and $T'$ are ($\leq 3$)-hypomorphic, then $T'_{\restriction S}\simeq \delta^+$ or $T'_{\restriction S}\simeq \delta^-$. We have $\{S\subseteq V, \ T'_{\restriction S} \simeq \delta^+ \}=\{S\subseteq V, \ T_{\restriction S} \simeq \delta^+ \}$, then $T'_{\restriction S}\simeq \delta^-$, $S\in U'^-$ and $U^-=U'^-$. So, for $X\subset V$, if $T_{\restriction X}$   is a diamond then $T'_{\restriction X} \simeq T_{\restriction X}$.\\
Now we prove that $T$ and $T'$ are $4$-hypomorphic. Let $X\subset V$ such that $|X|=4$. If $ T_{\restriction X} \simeq C_4$, then $c^{(3)}(T_{\restriction X})=2$. Since  $T$ and $T'$ are ($\leq 3$)-hypomorphic then $c^{(3)}(T'_{\restriction X})=2$, thus $T'_{\restriction X} \simeq T_{\restriction X} \simeq C_4$. The same, if $ T_{\restriction X} \simeq O_4$ then $T'_{\restriction X} \simeq T_{\restriction X} \simeq O_4$. So, $T'$ and $T$ are ($\leq 4$)-hypomorphic. Then, From Lemma \ref{hypomorphe}, $T'$ and $T$ are ($\leq 5$)-hypomorphic.
\endproof

From Claim \ref{3hyp4hyp}, it is sufficient to prove that $U^+=U'^+$.\\
For all $K\subseteq V$ with $\vert K\vert = k$, we have
$\{S\subseteq K : S\in U^+\}= D^+_4( T_{\restriction K})$ and
$\{S\subseteq K : S\in U'^+\}= D^+_4( T'_{\restriction K})$. \\
1)  Since $d^+_4(T_{\restriction K})=d^+_4(T'_{\restriction K})$ then
$\vert\{S\subseteq K : S\in U^+\}\vert=\vert\{S\subseteq K : S\in U'^+\}\vert$.
From Lemma \ref{particular mp}, we have $U^+=U'^+$.\\
2) We have $d^+_4(T_{\restriction K})\equiv d^+_4(T'_{\restriction K})$ (mod $p$)
for all $k$-element subsets $K$ of $V$, then $\vert\{S\subseteq K : S\in U^+\}\vert \equiv \vert\{S\subseteq K : S\in U'^+\}\vert$ (mod $p$).\\
a) Case 1. $p\geq5$, $t=4=[4]_p$, $k=[k_0,\dots]_p$ and $t_0=4\leq k_0$. From 1) of Theorem  \ref{thm js} we have $U^+=U'^+$.\\
Case 2. $p\geq5$, $t=4=[4]_p$, $k=[0,k_1,\dots]_p$. Since every tournament of cardinality  $\geq 5$ has at least a restriction of cardinality $4$ which is not a diamond, then from  2) of Theorem  \ref{thm js}, $U^+=U'^+$.\\
b)
Case 1. $p=3$, $t=4=[1,1]_p$, $k=[1,k_1,\dots,k_{k(p)}]_p$ and $t_1=1\leq k_1$. From 1) of Theorem  \ref{thm js} we have $U^+=U'^+$.\\
Case 2. $p=2$, $t=4=[0,0,1]_p$ and $k=[0,0,1,k_3,\dots,k_{k(p)}]_p$.\\
From 1) of Theorem  \ref{thm js} we have $U^+=U'^+$.\\
c)  We have
$p=2$, $t=4=[0,0,1]_p$, $k=[0,0,0,k_3,\dots,k_{k(p)}]_p$. Since every tournament of cardinality  $\geq 5$ has at least a restriction of cardinality $4$ which is not a diamond, and the fact that $T$ and $T'$ are  $3$-hypomorphic, then from  2) of Theorem  \ref{thm js},  $U^+=U'^+$,  thus $T'$ and $T$ are ($\leq 5$)-hypomorphic, or for all $4$-element subsets $S$ of V, $T_{\restriction S}$ is isomorphic to $\delta^+$ if and only if  $T'_{\restriction S}$ is isomorphic to $\delta^-$.
\endproof
\\

Given a digraph $S=(\{0,1,\dots ,m-1\},A)$, where $m\geq 1$
is an integer, for $i\in \{0,1,\dots ,m-1\}$ we associate a
digraph $G_{i}=(V_{i},A_{i})$, with $|V_i|\geq 1$, such
that the $V_{i}$'s are mutually disjoint. The
\textit{lexicographic sum} of $S$ by the digraphs $G_i$ or
simply the S-\textit{sum} of the $G_{i}$'s, is the digraph
denoted by $S(G_{0},G_{1},\dots ,G_{m-1})$ and defined on
the union of the $V_{i}$'s as follows: given $x\in V_{i}$
and $y\in V_{j}$, where $i,j\in \{0,1,\dots ,m-1\}$, $(x,y)$
is an arc of $S(G_{0},G_{1},\dots ,G_{m-1})$ if either $i=j$
and $(x,y)\in A_{i}$ or $i\neq j$ and $(i,j)\in A$: this
digraph replaces each vertex $i$ of $S$ by $G_{i}$.
We say that the vertex $i$ of $S$ is \textit{dilated by} $G_{i}$.\\

Let $h$ be a non-negative integer. The integers below are
considered modulo $2h+1$. The {\it circular tournament} $T_{2h+1}$  (see Figure $2$) is
defined on $\{0,1, \dots , 2h\}$ by :
${T_{2h+1}}_{\restriction \{0,1, \dots , h\}}$ is the usual total order on
$\{0,1, \dots , h\}, {T_{2h+1}}_{\restriction \{h + 1,\dots, 2h\}}$ is also
the usual order on $\{h + 1,h + 2, \dots , 2h\}$,  however
$ \{i + 1,i + 2, . . .\dots , h\}\longrightarrow_{T_{2h+1}}
i+h+1 \longrightarrow_{T_{2h+1}}\{0,1, \dots , i\}$ for every
$i\in \{0,1, \dots , h - 1\}$. A tournament $T$ is said to be an
element of $D(T_{2h+1})$ if $T$ is obtained by dilating each
vertex of $T_{2h+1}$ by a finite chain $p_{i}$, then
$T = T_{2h+1}(p_{0},p_{1},\dots ,p_{2h})$.
We recall that
$T_{2h+1}$ is indecomposable and $D(T_{2h+1})$ is the class
of finite tournaments without diamond  \cite{Lr}.\\

We define the tournament $\beta^+_6:= T_{3}(p_{0},p_{1},p_{2})$ with $p_0=(0<1<2)$, $p_1=(3<4)$ and $|p_2|=1$ (see Figure $3$). We set  $\beta^-_6:=(\beta^+_6)^*$.
For a tournament $T=(V,E)$, we set $B^+_6(T):=\{S\subseteq V : T_{\restriction S} \simeq \beta^+_6 \}$,
 $B^-_6(T):=\{S\subseteq V : T_{\restriction S}  \simeq \beta^-_6 \}$,
$b^+_6(T):=\mid B^+_6(T)\mid$   and $b^-_6(T):=\mid B^-_6(T)\mid$.\\

Two tournaments $T$ and $T'$ on the same vertex set $V$ are \textit{hereditarily isomorphic} if for all
$X\subseteq V$, $T_{\restriction X}$ and $T'_{\restriction X}$ are isomorphic \cite{BBN}.


%

\vspace{1.5cm}

\begin{figure}[h]
\setlength{\unitlength}{0.6cm}
\begin{center}
\begin{picture}(11,5)

 \put(6.0,7.0){$\bullet$\;0}   \put(8,6){$\bullet$\;1}
 \put(9,3.9){$\bullet$\;3}\put(9,5){$\bullet$\;2}
 \put(8,1){$\bullet$\;h-1}
 \put(6.1,0.1){$\bullet$\;h}
 \put(3.2,1){h+1\;$\bullet$}
  \put(2.1,1.9){h+2\;$\bullet$}
  \put(4,6.5){2h\;$\bullet$}
 \put(2.7,5){2h-1\;$\bullet$}
\put(6.2,7){\vector(2,-1){1.7}}
\put(6.2,7){\vector(3,-2){2.7}}
\put(6.2,7){\vector(1,-1){2.8}}
\put(6.2,7){\vector(1,-3){1.9}}
\put(6.2,7){\vector(0,-1){6.5}}
 \put(8.2,6){\vector(1,-1){0.8}}
 \put(8.2,6){\vector(1,-2){0.9}}
  \put(8.2,6){\vector(0,-1){4.6}}
   \put(8.1,6){\vector(-1,-3){1.8}}
   \put(8,6){\vector(-2,-3){3.1}}
\end{picture}
\end{center}
\vspace{-0.5cm}
\caption{Circular tournament $T_{2h+1}$}
\end{figure}
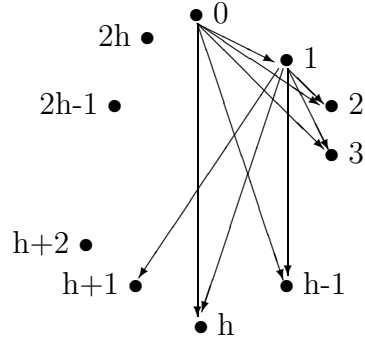

\vspace{0.33cm}

\vspace{-0.3cm}
\begin{figure}[h]
\begin{center}
\setlength{\unitlength}{0.4cm}
\begin{picture}(11,5)

\put(2.5,4.1){$(0<1<2)$}
\put(0.4,1){$5$} \put(9.1,1){$(3<4)$}

\put(8.8,1.1){\vector(-1,0){7.8}}

\put(1,1.3){\vector(3,2){3.9}}

\put(5.2,3.9){\vector(3,-2){3.9}}
\end{picture}
\vspace{-0.4cm}
\caption{$\beta_6^+$.}
\end{center}
\end{figure}
\vspace{- 0.5cm}

Let $G=(V,E)$  and $G'=(V,E')$ be two
$(\leq 2)$-hypomorphic digraphs. Denote  $D_{G,G'}$ the
binary relation on $V$ such that: for $x\in V$,
$x D_{G,G'}x$; and for $x \neq y\in V$, $x D_{G,G'}y$ if
there exists a sequence $x_{0} =x, . . . , x_{n} =y$ of
elements of $V$ satisfying $(x_{i}, x_{i+1})\in E$ if and
only if $(x_{i}, x_{i+1})\notin E'$, for all $i$,
$0 \leq i \leq n - 1$. The relation $D_{G,G'}$ is an
equivalence relation called {\it the difference relation},
its classes are called {\it difference classes}.

Using  difference classes, G. Lopez
\cite{L1,ls} showed that if $T$ and $T'$ are
($\leq 6$)-hypomorphic then $T$ and $T'$ are isomorphic.
One may deduce the next corollary.

\begin{corollary}\label{l1} (\cite{L1,ls}) Let $T$ and $T'$ be two tournaments. We have the following properties:\\
1)  If $T$ and $T'$ are $(\leq 6)$-hypomorphic then $T$ and $T'$ are hereditarily
 isomorphic.\\
2)  If for each
equivalence class $C$ of $D_{T,T'}$, $C$ is an interval of $T$ and $T'$, and
$T'_{\restriction C}$, $T_{\restriction C}$ are $(\leq 6)$-hypomorphic, then $T$ and $T'$ are hereditarily  isomorphic.
\end{corollary}

\begin{lemma}
\label{41}
\cite{L2}
Given two $(\leq 4)$-hypomorphic tournaments $T$
and $T'$, and  $C$  an equivalence class of $ D_{T,T'}$, then:\\
1) $C$ is an interval of $T'$ and $T$. \\
2)  Every $3$-cycle in
$T_{\restriction C}$ is reversed in $T'_{\restriction C}$.\\
3) There exists an integer $h\geq 0$ such that $T_{\restriction C}=T_{2h+1}(p_0,p_1,\dots,p_{2h})$ and $T'_{\restriction C}=T^*_{2h+1}(p'_0,p'_1,\dots,p'_{2h})$ with $p_i$, $p'_i$ are chains on the same basis, for all $i\in \{0,1,\dots ,2h\}$.
\end{lemma}

\begin{theorem}
Let $T$ and $T'$ be two  $(\leq 4)$-hypomorphic tournaments on the
same set $V$ of $v$ vertices. Let $p$ be a prime number and $k=[k_0,k_1,\dots ,k_{k(p)}]_p$  be an integer,  $6\leq k\leq v-6$.\\
1)  If $b^+_6(T_{\restriction K})=b^+_6(T'_{\restriction K})$ for all $k$-element subsets $K$ of $V$ then $T'$ and $T$ are hereditarily isomorphic.\\
2)   Assume  $b^+_6(T_{\restriction K})\equiv b^+_6(T'_{\restriction K})$ (mod $p$) for all $k$-element subsets $K$ of $V$.\\
a)   If $p\geq7$, and  $k_0\geq 6$ or $k_0=0$, then $T'$ and $T$ are hereditarily
 isomorphic.\\
b)   If ($p=5$, $k_0=1$ and $k_1\neq 0$) or ($p=3$, $k_0=0$ and $k_1=2$) or ($p=3$ and $k_0=k_1=0$)
or ($p=2$, $k_0=0$ and $k_1=k_2=1$), then $T'$ and $T$ are hereditarily
 isomorphic.
\end{theorem}

\Proof
Let $U^+:=\{S\subseteq V, \ T_{\restriction S} \simeq \beta^+_6 \}=B^+_6(T) $,
$U'^+:=B^+_6(T') $, $U^-:=\{S\subseteq V, \ T_{\restriction S} \simeq \beta^-_6 \}=B^-_6(T) $,
$U'^-:= B^-_6(T') $.\\
Every tournament of cardinality  $\geq 7$ has at least a restriction of cardinality $6$ which is not isomorphic to $\beta^+_6$ and $\beta^-_6$.
 Then for all cases,
similarly  to the proof of Theorem \ref{tournament}, we have $U^+=U'^+$.\\
Let $C$ be an equivalence class of $D_{T,T'}$, $S\in U^-$, $T_{\restriction S}\simeq \beta_6^-$. Since $T$ and $T'$ are ($\leq 3$)-hypomorphic, then $T'_{\restriction S}\simeq \beta^+_6$ or $T'_{\restriction S}\simeq \beta^-_6$. We have $\{S\subseteq V, \ T'_{\restriction S} \simeq \beta^+_6 \}=\{S\subseteq V, \ T_{\restriction S} \simeq \beta^+_6 \}$, then $T'_{\restriction S}\simeq \beta^-_6$, $S\in U'^-$ and $U^-=U'^-$. Let $X\subseteq C$ such that $|X|=6$;  if $T_X\simeq\beta^+_6$ then,  from 2) of Lemma \ref{41},  $T'_X\simeq\beta^-_6$, that is impossible, so $T_C$ and $T'_C$ has not a restriction of cardinality $6$ isomorphic to $\beta^+_6$ and $\beta^-_6$. \\
Now we will prove that $T_{\restriction C}$ and $T'_{\restriction C}$ are $(\leq 6)$-hypomorphic.\\
From 3) of Lemma \ref{41}, there exists an integer $h\geq 0$ such that
$T_{\restriction C}=T_{2h+1}(p_0,p_1,\dots,p_{2h})$, with $p_i$ is a chain and $a_i\in p_i$ for all $i\in \{0,1,\dots ,2h\}$ . Since $T_{\restriction C}$ hasn't a tournament isomorphic to $\beta_6^+$, then $h\leq3$. Indeed, if $h\geq4$, then $T_{\restriction\{a_0,a_1,a_2,a_3,a_4,a_{3+h}\}}\simeq  \beta_6^+$, and  $\{a_0,a_1,a_2\}$, $\{a_3,a_4\}$ are two intervals of  $T_{\restriction\{a_0,a_1,a_2,a_3,a_4,a_{3+h}\}}$, that is impossible.\\
a) If $h=3$, then $T_{\restriction C}=T_7$. Indeed, if $a_0,b_0\in V(p_0)$ then $T_{\restriction\{a_0,b_0,a_1,a_2,a_3,a_5\}}\simeq  \beta_6^+$, and $\{a_0,b_0,a_1\}$, $\{a_2,a_3\}$ are two intervals of  $T_{\restriction\{a_0,b_0,a_1,a_2,a_3,a_5\}}$, that is impossible.\\
b) If $h=2$, then $T_{\restriction C}=T_5$, or $T_{\restriction C}$ is obtained by dilating one vertex of $T_5$ by a chain of cardinality $2$. Indeed :

Case 1. $a_0,b_0,c_0\in V(p_0)$, then $T_{\restriction\{a_0,b_0,c_0,a_1,a_2,a_3\}}\simeq  \beta_6^+$ and $\{a_0,b_0,c_0\}$, $\{a_1,a_2\}$ are two intervals of  $T_{\restriction\{a_0,b_0,c_0,a_1,a_2,a_3\}}$, that is impossible.

Case 2. If  $a_i,b_i\in V(p_i)$ for all $i\in\{0,1\}$, then   $T_{\restriction\{a_0,b_0,a_1,b_1,a_3,a_4\}}\simeq  \beta_6^+$ and $\{a_0,b_0,a_4\}$, $\{a_1,b_1\}$ are two intervals of $T_{\restriction\{a_0,b_0,a_1,b_1,a_3,a_4\}}$, that is impossible.

Case 3. If  $a_i,b_i\in V(p_i)$ for all $i\in\{0,2\}$,  then $T_{\restriction\{a_0,b_0,a_1,a_2,b_2,a_4\}}\simeq  \beta_6^+$ and $\{a_0,b_0,a_1\}$, $\{a_2,b_2\}$ are two intervals of  $T_{\restriction\{a_0,b_0,a_1,a_2,b_2,a_4\}}$, that is impossible.\\
c) If $h=1$, then $T_{\restriction C}$ is obtained by dilating one vertex of $C_3$ by a chain or by dilating two or three vertices of $C_3$ by a chain of cardinality $2$.\\
d) If $h=0$, then $T_{\restriction C}$ is a chain.\\
In all cases, $T_{\restriction C}$ and $T'_{\restriction C}$ are ($\leq 6$)-hypomorphic. From 1) of Lemma \ref{41}, $C$ is an interval of $T'$ and $T$. Then, from 2) of Corollary \ref{l1},  $T$ and $T'$ are hereditarily isomorphic.
\endproof


\end{document}